\documentclass[12pt,a4paper]{article}
\usepackage[utf8]{inputenc}
\usepackage{amsmath}
\usepackage{amsfonts}
\usepackage{amssymb}
\usepackage{graphicx}
\usepackage[left=2.5cm,right=2.5cm,top=2.5cm,bottom=2.5cm]{geometry}

\usepackage{authblk}
\usepackage{array}
\usepackage{lineno}
\usepackage{hyperref}

\author[1]{Lauren D Smith}
\affil[1]{Department of Mathematics, The University of Auckland, Auckland 1142, New Zealand}
\date{}
\title{Optimal illness policy for an unethical daycare center}

\graphicspath{{figures/}}

\begin{document}

\maketitle

\begin{abstract}
While businesses are typically more profitable if their workers and communities are minimally exposed to diseases, the same is not true for daycare centers. Here it is shown that a daycare center could maximize its profits by maintaining a population of sick children within the center, with the intention to infect more children who then do not attend. Through a modification of the Susceptible-Infected-Recovered (SIR) model for disease spread we find the optimal number of sick children who should be kept within the center to maximize profits. We show that as disease infectiousness increases, the optimal attendance rate of sick children approaches zero, while the potential profit increases. 
\end{abstract}

\section{Introduction}

A huge number of children and infants attend daycare centers. We as a community trust that daycare centers will act ethically, putting the needs and health of children above all else. While this is likely the case in most centers, in this article it will be shown that daycare centers have a financial incentive to spread diseases among children and infants, which is a clear conflict between ethics and profit.

Providing paid sick leave is widely regarded to yield a net profit for businesses \cite{asfaw_potential_2017, lohaus_presenteeism_2019, vander_weerdt_is_2023}, i.e., the loss in productivity from illnesses spreading among the workforce is greater than the cost of paid sick leave. Businesses are therefore financially incentivized to minimize disease spread, which is beneficial for society as a whole.
 In this article it will be shown that daycare centers are financially incentivized to keep ill children in the center with the purpose of infecting other children, and that there is an optimal proportion of sick children that should attend to maximize center profits. The incentive exists under the following assumptions (which are true at numerous daycare centers globally as identified from policy documents and ``parent handbooks''):
\begin{enumerate}
	\item Teachers on casual work contracts are sent home (unpaid) when child attendance is low enough to do so (retaining only the minimum number of teachers for government imposed teacher:child ratios).
	\item Children pay regardless of whether they attend or not (payment by enrollment rather than attendance).
	
	\item Teachers do not receive paid sick leave.
\end{enumerate}
The first assumption is the most important. By maintaining only the minimum number of teachers, the daycare center can minimize its staffing costs by minimizing the attendance of its students (maximizing absence). If children pay regardless of attendance (Assumption 2), this further increases the potential profit from student absence. Paying regardless of attendance provides financial stability for daycare centers, but also provides a financial incentive for parents to send sick children to daycare (otherwise parents may not be paid). 
 By not providing paid sick leave to teachers (Assumption 3) there is no financial incentive to minimize staff illness. While not included in the modeling here, Assumption 3 also provides a financial incentive for sick teachers to continue working, which will further spread illness. Essentially, daycare centers with the above assumptions can maximize profits by maximizing the number of days that students are absent from the center, which can be achieved through their illness policy (both the policy itself and selective enforcement). At one extreme, in a daycare center with zero tolerance for illness, there will be very little spread of disease, and,  hence, a low number of days absent. At the other extreme, a daycare center that allows all ill children to attend will also have a small number of days absent, because ill children are never kept home. Between these two extremes is an optimal attendance rate at which the disease is still able to spread and some of the ill children are sent home, with a maximal total absence time across all children.
 
 Utilizing a modification of the well-known \textit{Susceptible-Infected-Recovered} (SIR) compartment model for infection spread \cite{kermack_contribution_1927, hethcote_mathematics_2000}, we determine the optimal attendance rate to maximize profits, and calculate the potential savings that can be made by employing the optimal illness policy. It is found that more savings can be made from highly infectious diseases such as measles compared to less infectious diseases. It is also found that the optimal attendance rate tends toward zero as the infectiousness of diseases increases, which is a positive for a society that is trying to minimize disease.

While there is no evidence that daycare centers are deliberately keeping sick children in classes to spread diseases, the practice of sending staff home when child attendance is low (Assumption 1) has been previously reported \cite{rush_abc_case_2006}, and is explicitly written in some daycare center policy documents (typically as part of ``late arrival'' or ``drop-off' policies'). There are many more daycare centers that have the same ``late-arrival'' or ``drop-off'' policies without explicitly stating that staff will be sent home if attendance is low. The practice is also widely reported on internet forums by early childhood education (ECE) teachers. It has also been reported that daycare centers over-enroll children (i.e., under-staff teachers) with the assumption that a proportion of students will not attend due to illness \cite{oece_2023}.  
All daycare centers that either reactively send staff home or preemptively under-staff are taking advantage of the potential cost-savings proposed here, but are unlikely to be doing so optimally.

\section{The model}

\begin{figure}[tbp]
\centering
\includegraphics[width=0.8\columnwidth]{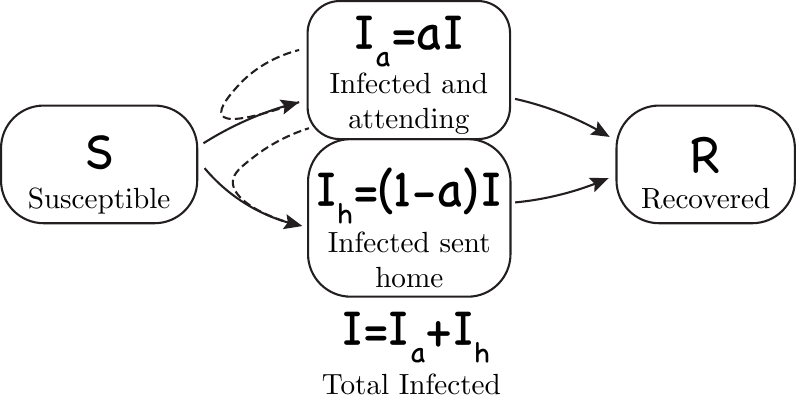}
\caption{
Schematic of the modified SIR compartment model. The total population is divided into susceptible ($S$), infected ($I$) and recovered ($R$) individuals, and the infected population is further divided into those who attend daycare ($I_a$) and those who are sent home ($I_h$). Susceptible individuals become infected due to interactions with individuals who are infected and attending. Infected individuals gradually recover.
}
\label{fig:SIR_schematic}
\end{figure}

We consider the scenario in which a new illness is introduced to a daycare center, such that the majority of children are initially susceptible to the illness. 
In the standard SIR model, the total population ($N$) of children is divided into three categories: those who are susceptible and have no immunity to the disease ($S$), those who are currently infected and can infect others ($I$), and those who have recovered from the disease and have gained immunity ($R$). Children are able to move between categories as time progresses due to becoming infected or recovering from the disease. Here we modify the standard model such that the infected population is divided into two groups; those who continue to attend the daycare (and who can then infect other children) and those who stay home. We define $a$ as being the proportion of infected children who attend daycare. The value of $a$ can be influenced by a number of factors, including the actions of parents and their willingness to send sick children to daycare, as well as how quickly and frequently a daycare will send sick children home or not accept sick children into their care. A daycare with a very strict illness policy, that tries to minimize the impact of a disease on the greater community, would have $a\approx 0$, whereas a daycare with a very lenient illness policy could have $a \approx 1$. We denote by $I_a = a I$ the number of infected children that attend daycare, and denote by $I_h = (1-a) I$ the number of infected children who stay home.
Fig.~\ref{fig:SIR_schematic} shows a schematic diagram of the situation, such that susceptible children ($S$) become infected due to interactions with attending infected children ($I_a$). The newly infected children are split (based on $a$) into those who attend daycare ($I_a$) and those who stay at home ($I_h$). All infected individuals gradually recover from the illness, moving them into the recovered ($R$) category, at which point they are assumed to be immune to the disease.
 Mathematically, the dynamics can be described by the differential equations
 \begin{linenomath}
\begin{align}
\frac{dS}{dt} &= - \frac{\beta}{N} I_a \,S = - a \frac{\beta}{N} I \, S,  \nonumber \\
\frac{dI}{dt} &= \frac{\beta}{N} I_a \, S - \gamma I = a \frac{\beta}{N} I \, S - \gamma I, \label{eq:SIR} \\
\frac{dR}{dt} &= \gamma I, \nonumber
\end{align}
\end{linenomath}
where $\beta$ is the infection rate of the disease and $\gamma$ is the recovery rate of the disease. We see that $a$ directly scales the infectiousness of the disease $\beta$. With $a=1$, the system (\ref{eq:SIR}) is the standard SIR model. 

The model (\ref{eq:SIR}) assumes a \textit{well-mixed} population, in the sense that all individuals interact with each other all the time, which can be reasonably expected in a classroom with young children and infants. We have also assumed that teachers have a negligible influence on the infection of children, which we propose based on assumption that teachers are more hygienic than children, meaning they are less prone to spreading infections. There are also fewer teachers than children within a daycare center, meaning fewer opportunities to spread a disease from a teacher to a child.

The typical dynamics of the model (\ref{eq:SIR}) are illustrated in Fig.~\ref{fig:sims_and_Th}(a-d) for a range of attendance rates $a$. The disease outbreak starts with a single initially infected child ($I(0)=1$) among a group of $N=100$ children. Note that all populations approach limiting values as time progresses, and not all children become infected over the course of the disease. From a public health point-of-view, the optimal attendance rate $a$ corresponds to minimizing the total number of infections over the course of the disease, which means a value of $a$ as close to zero as possible. For example, the final value of $S$ in Fig.~\ref{fig:sims_and_Th}(a) with $a=0.3$ is much greater than that in Fig.~\ref{fig:sims_and_Th}(d) with $a=0.8$, meaning less total infections for $a=0.3$. However, as we will demonstrate, minimizing infections is not the most profitable scenario for a daycare center.

\section{Maximizing daycare center profits}

\begin{figure}[tbp]
\centering
\includegraphics[width=\columnwidth]{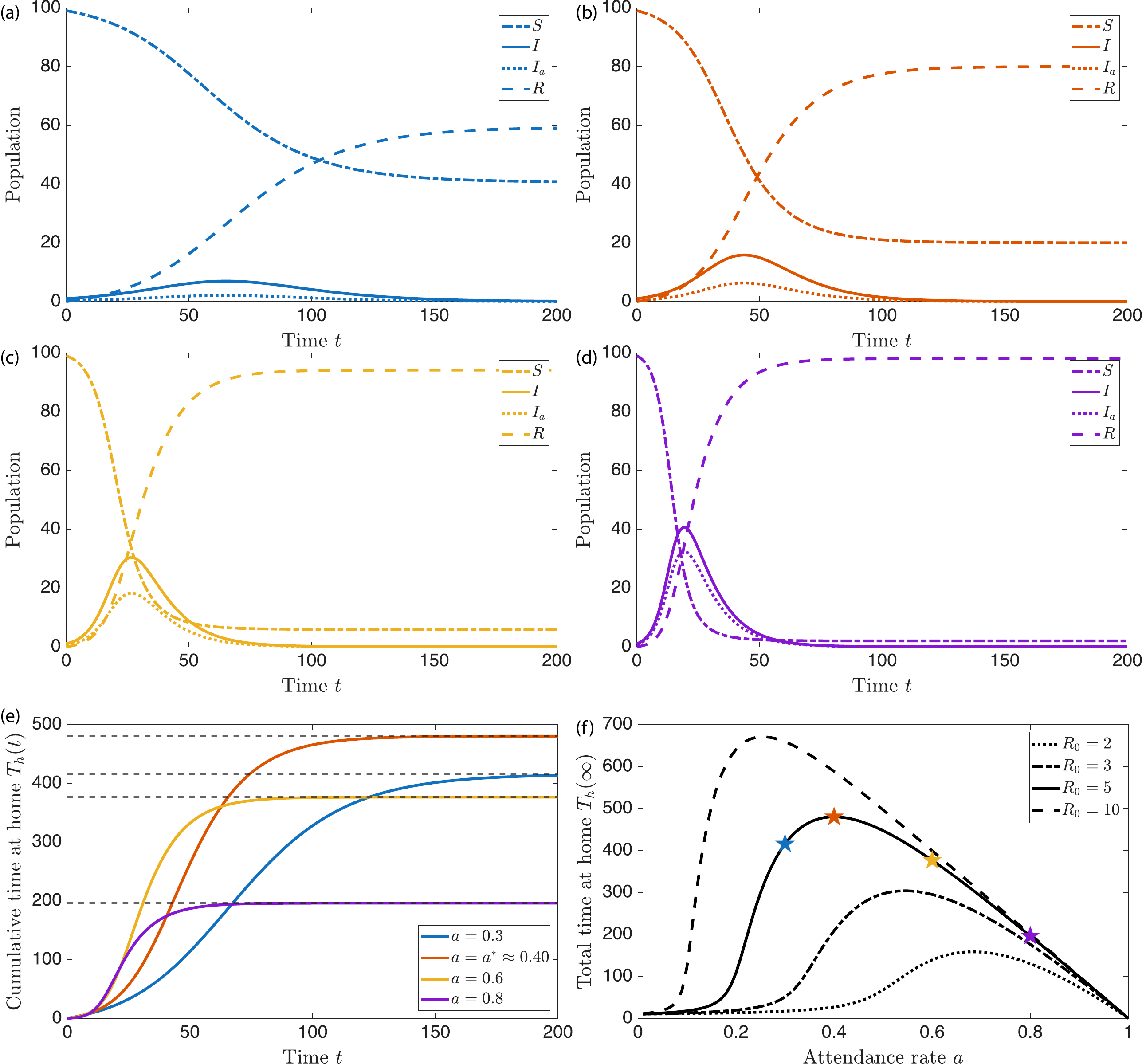}
\caption{
(a-d) Time evolution of the components $S$ (dot-dashed), $I$ (solid), $I_a$ (dotted) and $R$ (dashed) of the modified SIR model (\ref{eq:SIR}) with $N=100$ children, infectivity $\beta = 0.5$, recovery rate $\gamma = 0.1$, basic reproduction number $R_0 = \beta/\gamma = 5$, initial populations $S(0) = 99$, $I(0) = 1$, $R(0) = 0$, and varied attendance rate $a$. (a)~$a=0.3$, (b)~$a=a^*\approx 0.40$, (c)~$a=0.6$, (d)~$a=0.8$. (e)~The cumulative time that children are kept home $T_h(t)$ (\ref{eq:Th_t}) for the same parameters as in (a-d), with corresponding colors. The dashed horizontal lines show respective long-term theoretical values $T_h(\infty)$ based on (\ref{eq:Th_infty}) and (\ref{eq:r_infty}). (f)~Dependence  of the total time at home over the course of the disease $T_h(\infty)$ on the attendance parameter $a$ for a range of $R_0$ values (all other parameters the same as in (a-e)). The cases directly corresponding to the parameters in (a-e) are highlighted with stars in the respective colors from (a-d).
}
\label{fig:sims_and_Th}
\end{figure}

A daycare center that is maximizing profits would minimize its staffing costs by having the smallest number of teachers required to meet the required staff to child ratio at any given time. Therefore, profits are maximized by maximizing the total amount of time that children are kept home over the course of the disease.
We let $T_h(t)$ denote the cumulative amount of time that children have been kept home from the start of the outbreak up until the time $t$, which can be calculated as
\begin{equation} \label{eq:Th_t}
T_h(t) = \int_0^t I_h(u) du =(1-a) \int_0^t I(u) du,
\end{equation}
and is shown in Fig.~\ref{fig:sims_and_Th}(e) for a range of $a$ values (as in Fig.~\ref{fig:sims_and_Th}(a-d)) and all other parameters kept fixed. Like the populations $S$, $I$, and $R$, the cumulative time at home $T_h(t)$ approaches a limiting value, $T_h(\infty) = \lim_{t\to\infty} T_h(t)$, which is the quantity that should be maximized to maximize daycare profits. Note that of the four cases shown in Fig.~\ref{fig:sims_and_Th}(a-d), the case $a=0.8$ (Fig.~\ref{fig:sims_and_Th}(d)) has the largest total number of children that are infected over the course of the disease (98\% infected), but the lowest long-term time children are kept home ($T_h(\infty) = 196$). Comparatively, the case with $a=a^*\approx 0.40$ (Fig.~\ref{fig:sims_and_Th}(b)) which has maximal $T_h(\infty) = 480$, has fewer total infections (80\% infected).

The long-term time at home $T_h(\infty)$ is proportional to the fraction of infected children who are kept home $(1-a)$, the total number of children who are infected over the course of the disease, and the typical length of the infection $1/\gamma$. Since every infected individual eventually recovers, the total number of children who are infected over the course of the disease is equal to the long-term population of recovered individuals, $R_\infty = \lim_{t\to \infty} R(t)$. Therefore, we find that
\begin{equation} \label{eq:Th_infty}
T_h(\infty) = \frac{(1-a) R_\infty}{\gamma}.
\end{equation}
The value $R_\infty$ depends on $a$ through the dynamics (\ref{eq:SIR}), and can be written as $R_\infty = N r_\infty$, where $r_\infty$ is the proportion of the population that is recovered in the long term. 
Borrowing from the standard SIR model, we get an equation for $r_\infty$ \cite{wang_application_2010, kroger_analytical_2020},
\begin{equation} \label{eq:r_infty}
r_\infty = 1 + \frac{1}{a R_0} W \left( -s_0 a R_0 e^{- a R_0 (1-r_0)} \right),
\end{equation}
where $R_0 = \beta/\gamma$ is the \textit{basic reproduction number} of the disease, $s_0 = S(0)/N$ is the initial proportion of susceptible individuals, $r_0 = R(0)/N$ is the initial proportion of recovered individuals, and $W(x)$ is the Lambert $W$ function. We will assume that initially 1\% of the child population is infected, and no children have immunity, meaning $s_0 = 0.99$ and $r_0 = 0$. For a given disease, which has inherent infectiousness $\beta$, recovery time $\gamma$, and $R_0 =\beta/\gamma$, the expressions (\ref{eq:Th_infty}) and (\ref{eq:r_infty}) depend only on the attendance rate $a$. Fig.~\ref{fig:sims_and_Th}(f) shows $T_h(\infty)$ as a function of $a$ for a range of $R_0$ values (with $\gamma = 0.1$ kept fixed). The solid curve with $R_0=5$ corresponds to the parameters in Fig.~\ref{fig:sims_and_Th}(a-e), and the colored stars correspond to the simulations of (\ref{eq:SIR}) in Fig.~\ref{fig:sims_and_Th}(a-d). 

\begin{figure}[tbp]
\centering
\includegraphics[width=0.7\columnwidth]{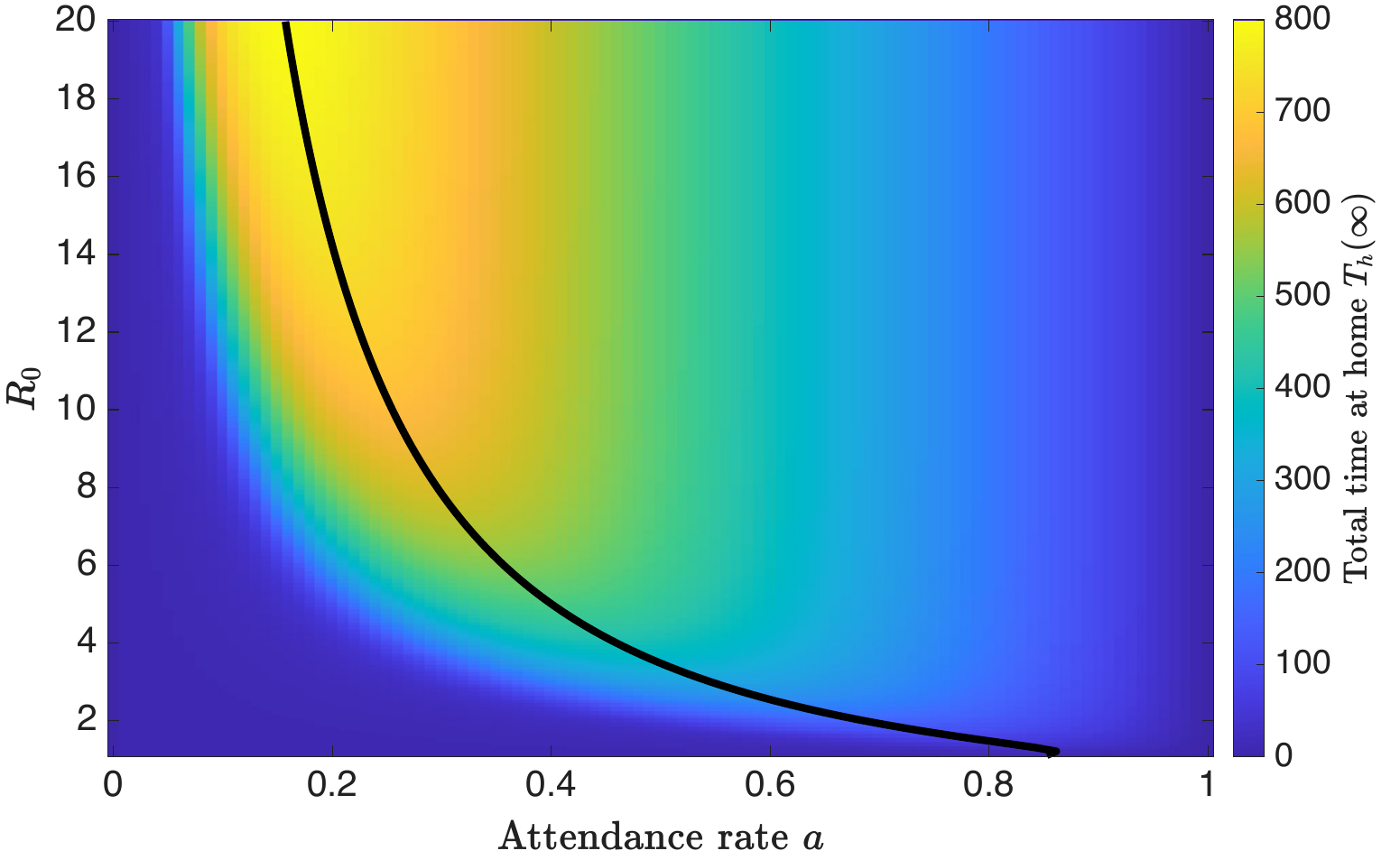}
\caption{
Long-term time at home $T_h(\infty)$ across the $(a,R_0)$ parameter space, with $\gamma= 0.1$ and $N=100$ kept fixed. The black curve shows the attendance rates $a^*$ with maximal $T_h(\infty)$ for fixed $R_0$.
}
\label{fig:astar_Th}
\end{figure}

Maximizing $T_h(\infty)$ (and, hence, daycare profits) amounts to solving the simple numerical root-finding problem 
\begin{equation}
\frac{d T_h(\infty)}{da} = 0,
\end{equation}
to find the optimal attendance rate $a^*$. The dependence of $a^*$ on the basic reproduction number $R_0$ of the disease is shown by the black curve in Fig.~\ref{fig:astar_Th}, which also shows the long-term time at home $T_h(\infty)$ as a function of $a$ and $R_0$. Fig.~\ref{fig:astar_Th} shows that as $R_0$ increases, the maximum time at home $T_h(\infty)$ increases, and the optimal attendance rate $a^*$ at which the maximum occurs decreases towards zero. The trend towards zero of $a^*$ as $R_0$ increases means that daycare profitability more strongly aligns with public health outcomes for more infectious diseases.

Table~\ref{tab:diseases} shows the optimal attendance rates $a^*$ for a variety of diseases. Also shown is the potential savings (in worked days) that can be made in staffing by adopting the optimal attendance rate, assuming a teacher to child ratio of 1:6 and $N=100$ children. The savings are calculated as
\begin{equation}
\text{Savings} = T_h(\infty)|_{a^*} \times \frac{1}{6} \times \frac{5}{7},
\end{equation}
which accounts for the assumed teacher:child ratio and that children only attend daycare on 5 out of 7 days per week. As discussed previously, diseases such as measles with higher $R_0$ values have higher maximum $T_h(\infty)$ values, which in turn offer the highest potential savings.

\begin{table}
\begin{center}
\begin{tabular}{|c|c|c|c|c|}
\hline 
Disease & $R_0$ & $a^*$ & $T_h(\infty)|_{a^*}$ & Staff savings (days) \\ 
\hline 
\hline
Measles & 12-18 \cite{guerra_basic_2017}  & 0.19 & 753 & 90 \\ 
\hline 
COVID-19 (Omicron) & 9.5 \cite{liu_effective_2022} & 0.26 & 659 & 78 \\ 
\hline 
COVID-19 (ancestral) & 2.4-3.4 \cite{billah_reproductive_2020} & 0.56 & 291 & 35 \\ 
\hline 
Influenza (seasonal) & 1.2-1.4 \cite{chowell_seasonal_2008} & 0.84 & 36.7 & 4.4 \\ 
\hline 
\end{tabular} 
\caption{
For a range of diseases: optimal attendance rates $a^*$, the corresponding maximum total time that children are kept home over the course of the disease $T_h(\infty)|_{a^*}$, and the potential staffing reductions that can be made. Results are for $N=100$ children, an assumed teacher to child ratio of 1:6, and a recovery rate $\gamma = 0.1$ (i.e., typical recovery in $1/\gamma = 10$ days).
}
\end{center}
\label{tab:diseases}
\end{table}

Note that for a typical business, they would want to maximize productivity, which means minimizing the total time that workers spend at home $T_h(\infty)$ while recovering. From Fig.~\ref{fig:sims_and_Th}(f) and Fig.~\ref{fig:astar_Th} we see that $T_h(\infty)$ is minimized (and equal to zero) with either $a=0$ (no attendance of ill workers) or $a=1$ (full attendance of ill workers). But since ill workers are generally less productive than healthy workers, it is best for the business to have $a=0$. Paid sick leave provides an incentive for ill workers to stay home, bringing $a$ closer to zero, and it has been shown that the increase in productivity gained by ill workers staying home outweighs the cost of paid sick leave \cite{asfaw_potential_2017, lohaus_presenteeism_2019, vander_weerdt_is_2023}.

\section{Discussion}

In summary, we have shown that under reasonable assumptions daycare centers are financially incentivized to keep a non-zero proportion of sick children in their care so as to spread disease to other children. We have found the optimal attendance rate, which depends on the \textit{basic reproduction number} $R_0$ of the disease, and the potential cost savings that could be made in terms of reduced staff hours.

The SIR model is particularly simple, which is good for highlighting the general trends presented in this article. If more accurate results are desired one could implement a more detailed model which incorporates more complicated dynamics such as ``exposed'' populations that are not yet infectious (the SEIR model) \cite{hethcote_mathematics_2000}, as well as stochasticity \cite{tornatore_stability_2005, ji_threshold_2014}, or a model which incorporates the discreteness of the true situation, e.g., as a branching process \cite{hendy_mathematical_2021}.

The model here assumes that teachers do not receive paid sick leave (Assumption 3). We expect that the financial incentive for sick students to attend daycare centers will persist even when paid sick leave is included. There would be an initial penalty for paying the sick leave, but once all sick leave is exhausted the penalty is removed.

This article is not intended as a recipe for unethical daycare centers, but rather to illustrate the financial incentive that exists for daycare centers to propagate diseases among children, which would lead to more infections of at-risk populations in the wider community. The financial incentive can be reduced or removed by addressing the three key assumptions presented at the beginning of this article. In particular, removing the practice of sending teachers home to meet minimum required teacher:child ratios would eliminate the financial incentive entirely. In addition, paid sick-leave for teachers would provide a financial disincentive for fostering disease within centers.

\vspace{12pt}

\noindent \textbf{Data Availability:} The MATLAB code required to reproduce the results presented is available in the GitHub repository \url{https://github.com/LaurenSmithMaths/UnethicalDaycare}

%\bibliographystyle{ieeetr} %%%% .BST file
%
%\bibliography{bibliography.bib} %%%%% .Bib file

\begin{thebibliography}{10}

\bibitem{asfaw_potential_2017}
A.~Asfaw, R.~Rosa, and R.~Pana-Cryan, ``Potential {Economic} {Benefits} of
  {Paid} {Sick} {Leave} in {Reducing} {Absenteeism} {Related} to the {Spread}
  of {Influenza}-{Like} {Illness},'' {\em Journal of occupational and
  environmental medicine}, vol.~59, pp.~822--829, Sept. 2017.

\bibitem{lohaus_presenteeism_2019}
D.~Lohaus and W.~Habermann, ``Presenteeism: {A} review and research
  directions,'' {\em Human Resource Management Review}, vol.~29, pp.~43--58,
  Mar. 2019.

\bibitem{vander_weerdt_is_2023}
C.~Vander~Weerdt, P.~Stoddard-Dare, and L.~DeRigne, ``Is paid sick leave bad
  for business? {A} systematic review,'' {\em American Journal of Industrial
  Medicine}, vol.~66, no.~6, pp.~429--440, 2023.
\newblock \_eprint: https://onlinelibrary.wiley.com/doi/pdf/10.1002/ajim.23469.

\bibitem{kermack_contribution_1927}
W.~O. Kermack and A.~G. McKendrick, ``A contribution to the mathematical theory
  of epidemics,'' {\em Proceedings of the Royal Society of London. Series A,
  Containing Papers of a Mathematical and Physical Character}, vol.~115,
  pp.~700--721, Jan. 1927.
\newblock Publisher: Royal Society.

\bibitem{hethcote_mathematics_2000}
H.~W. Hethcote, ``The {Mathematics} of {Infectious} {Diseases},'' {\em SIAM
  Review}, vol.~42, pp.~599--653, Jan. 2000.
\newblock Publisher: Society for Industrial and Applied Mathematics.

\bibitem{rush_abc_case_2006}
E.~Rush and C.~Downie, ``{ABC Learning Centres} - {A} case study of
  {A}ustralia's largest child care corporation.'' Report of the Australia
  Institute, June 2006.

\bibitem{oece_2023}
N.~Z. Office~of Early Childhood~Education, ``The quality of early childhood
  education provided to children - 3-yearly teacher survey results,'' 2023.

\bibitem{wang_application_2010}
F.~Wang, ``Application of the {Lambert} {W} {Function} to the {SIR} {Epidemic}
  {Model},'' {\em The College Mathematics Journal}, vol.~41, pp.~156--159, Mar.
  2010.
\newblock Publisher: Taylor \& Francis \_eprint:
  https://doi.org/10.4169/074683410X480276.

\bibitem{kroger_analytical_2020}
M.~Kr{\"o}ger and R.~Schlickeiser, ``Analytical solution of the {SIR}-model for
  the temporal evolution of epidemics. {Part} {A}: time-independent
  reproduction factor,'' {\em Journal of Physics A: Mathematical and
  Theoretical}, vol.~53, p.~505601, Nov. 2020.
\newblock Publisher: IOP Publishing.

\bibitem{guerra_basic_2017}
F.~M. Guerra, S.~Bolotin, G.~Lim, J.~Heffernan, S.~L. Deeks, Y.~Li, and N.~S.
  Crowcroft, ``The basic reproduction number ({R0}) of measles: a systematic
  review,'' {\em The Lancet Infectious Diseases}, vol.~17, pp.~e420--e428, Dec.
  2017.
\newblock Publisher: Elsevier.

\bibitem{liu_effective_2022}
Y.~Liu and J.~Rockl{\"o}v, ``The effective reproductive number of the {Omicron}
  variant of {SARS}-{CoV}-2 is several times relative to {Delta},'' {\em
  Journal of Travel Medicine}, vol.~29, p.~taac037, Mar. 2022.

\bibitem{billah_reproductive_2020}
M.~A. Billah, M.~M. Miah, and M.~N. Khan, ``Reproductive number of coronavirus:
  {A} systematic review and meta-analysis based on global level evidence,''
  {\em PLoS ONE}, vol.~15, p.~e0242128, Nov. 2020.

\bibitem{chowell_seasonal_2008}
G.~Chowell, M.~A. Miller, and C.~Viboud, ``Seasonal influenza in the {United}
  {States}, {France}, and {Australia}: transmission and prospects for
  control,'' {\em Epidemiology and Infection}, vol.~136, pp.~852--864, June
  2008.

\bibitem{tornatore_stability_2005}
E.~Tornatore, S.~Maria~Buccellato, and P.~Vetro, ``Stability of a stochastic
  {SIR} system,'' {\em Physica A: Statistical Mechanics and its Applications},
  vol.~354, pp.~111--126, Aug. 2005.

\bibitem{ji_threshold_2014}
C.~Ji and D.~Jiang, ``Threshold behaviour of a stochastic {SIR} model,'' {\em
  Applied Mathematical Modelling}, vol.~38, pp.~5067--5079, Nov. 2014.

\bibitem{hendy_mathematical_2021}
S.~Hendy, N.~Steyn, A.~James, M.~J. Plank, K.~Hannah, R.~N. Binny, and
  A.~Lustig, ``Mathematical modelling to inform {New} {Zealand}'s {COVID}-19
  response,'' {\em Journal of the Royal Society of New Zealand}, vol.~51,
  pp.~S86--S106, May 2021.
\newblock Publisher: Taylor \& Francis \_eprint:
  https://doi.org/10.1080/03036758.2021.1876111.

\end{thebibliography}

\end{document}